\documentclass[12pt]{article}%
\usepackage{amsfonts}
\usepackage{sw20bams}
\usepackage{amsmath}
\usepackage{amssymb}
\usepackage{graphicx}%
\providecommand{\U}[1]{\protect\rule{.1in}{.1in}}
\begin{document}

\title{A translation of Henri Joris' \textquotedblleft Le chasseur perdu dans la
for\^{e}t\textquotedblright\ (1980)}
\author{Steven Finch}
\date{October 1, 2019}
\maketitle

\begin{abstract}
This is an English translation of Henri Joris' article \textquotedblleft Le
chasseur perdu dans la for\^{e}t (Un probl\`{e}me de g\'{e}om\'{e}trie
plane)\textquotedblright\ that appeared in \textit{Elemente der Mathematik} v.
35 (1980) n. 1, 1--14. \ Given a point $P$ and a line $L$ in the plane, what
is the shortest search path to find $L$, given its distance but not its
direction from $P$? The shortest search path was described by Isbell (1957),
but a complete and detailed proof was not published until Joris (1980).
\ I\ am thankful to Natalya Pluzhnikov for her dedicated work and to the Swiss
Mathematical Society for permission to post this translation on the arXiv.

\end{abstract}

1. A hunter wandered into the woods and became disoriented. After awhile, he
found a sign with information about a road that passes exactly $1$ kilometer
away. Unfortunately, the tree on which the sign was attached had fallen to the
ground, and the hunter possessed no idea in what direction he should travel to
reach the road. He decided to walk $1$ kilometer straight in an arbitrary
direction and then continue along a circle with center at the sign location.
Therefore, he is certain to find the road after walking at most $1+2\pi$ kilometers.

This situation is described in a problem presented at a mathematical
competition for high school students, which are so popular in the United
States. The question is: \textit{Suppose the road is straight. What is the
shortest search path that the hunter would have taken, had he thought about
everything a little more carefully?} \ The answer probably is:
\textquotedblleft Instead of traversing the full circumference of the circle,
he would replace the last quarter by a straight path of length $1$ km tangent
to the circle (Figure 1), saving $\pi/2-1\approx0.57$ km\textquotedblright.%
\begin{figure}[ptb]%
\centering
\includegraphics[
height=2.5786in,
width=6.1668in
]%
{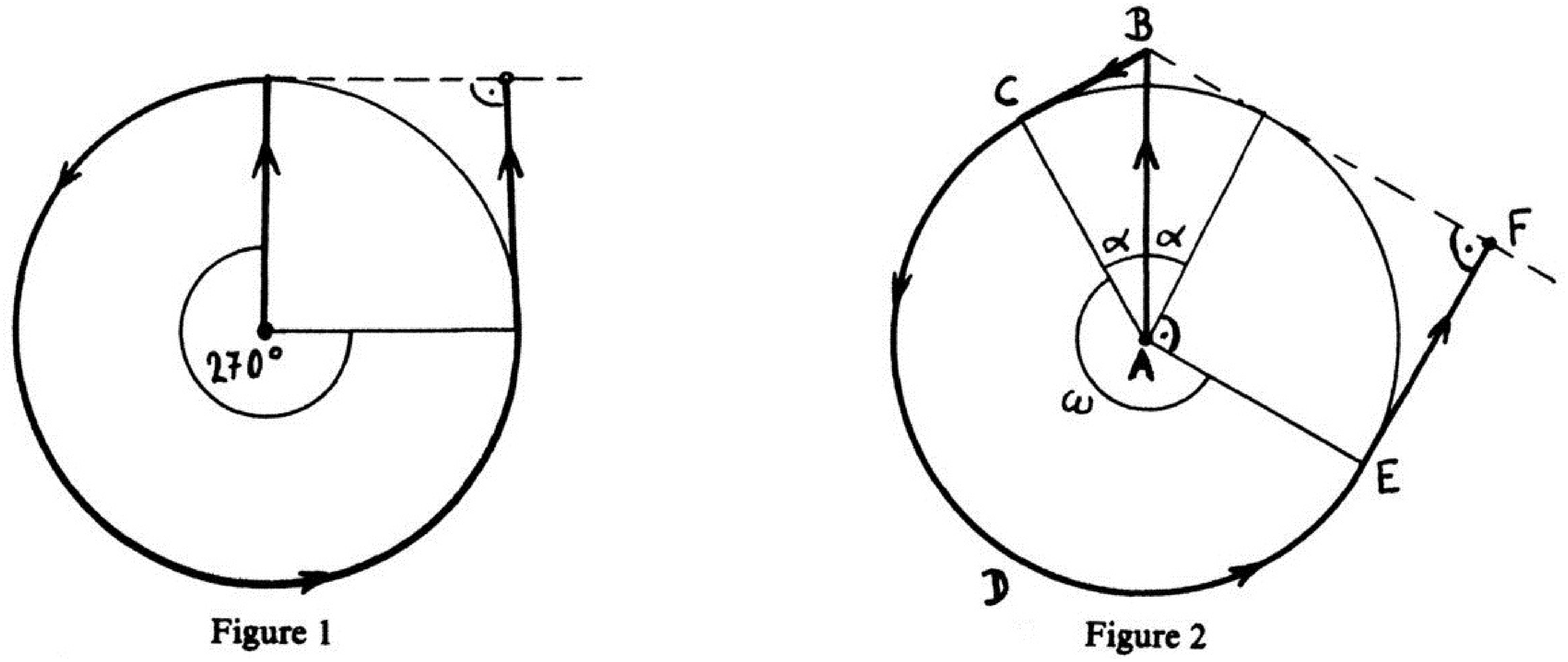}%
\end{figure}

However, he can make his way even shorter if he first continues to walk in the
direction of the radius beyond the circle and then return to the circle
tangentially (Figure 2). This will give the path $ABCDEF$ of length%

\[
\frac{1}{\cos\alpha}+\tan\alpha+\omega+1=\frac{1}{\cos\alpha}+\tan\alpha
+\frac{3}{2}\pi-2\alpha+1=\ell(\alpha).
\]
$\ell(\alpha)$ attains its minimum at $\alpha=30^{\circ}=\pi/6$,%

\begin{equation}
\ell\left(  \frac{\pi}{6}\right)  =\frac{7}{6}\pi+1+\sqrt{3}.
\end{equation}

2. If we disassociate from the hunter and the forest, and rephrase in terms of
plane geometry, we come to the following problem (P):

\textit{Given a circle of radius 1, find the shortest path that starts at the
center of the circle and intersects all the tangents to the circle.}

The following theorem holds:

\bigskip

\textbf{Theorem 1}. \textit{A solution of (P) is the path of Figure 2 with
}$\alpha=\pi/6$\textit{. Any other solution is obtained from it by rotations
and reflections of the plane that leave the circle in its place.}

\bigskip

If we remove from (P) the condition about the starting point, we obtain

\bigskip

\textbf{Theorem 2}. \textit{The shortest path that intersects all the tangents
of a given circle consists of a semicircle (subset of the given circle)
extended on each side by a segment of the tangent of length of the radius}
(Figure 3).

\bigskip

(If the radius is of length $1$, the length of the path is $\pi+2$.)%
\begin{figure}[ptb]%
\centering
\includegraphics[
height=2.6999in,
width=3.2818in
]%
{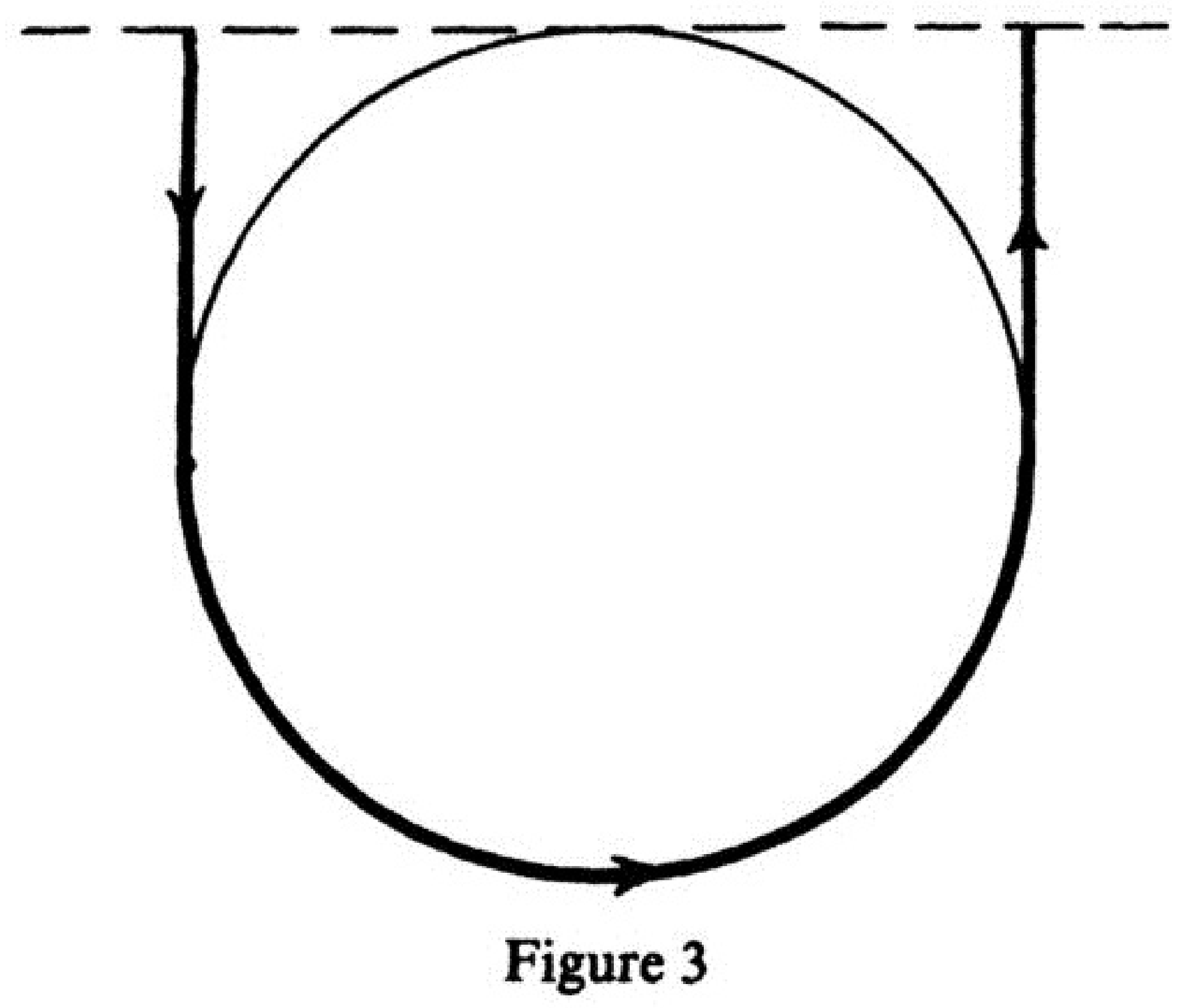}%
\end{figure}

Although these theorems, especially the second, look quite elementary and
plausible, I could not find simple and totally elementary proofs of them. This
is due to the fact that one must take into consideration all continuous and
rectifiable curves that satisfy the conditions, and also that the solutions
must be curves composed of straight line segments and circular arcs.

\bigskip

3. In what follows, we shall prove these two theorems. The idea is to show
that each minimal curve is one of those indicated in the theorems. First we
prove the existence of minimal curves.

Let $m$ be the infimum of the lengths of the curves that meet all the
tangents, and $\mathfrak{C}_{n}$ a minimal sequence of such curves; then
$\ell(\mathfrak{C}_{n})=($length of $\mathfrak{C}_{n})\rightarrow m$. Each
$\mathfrak{C}_{n}$ can be parametrized by $f_{n}:[0,1]\rightarrow E_{2}$, the
parameter being proportional to the arc length. The $\ell(\mathfrak{C}_{n})$
are bounded; therefore, all the $\mathfrak{C}_{n}$ are contained in a closed
square, and%
\[
\operatorname{distance}\left(  f_{n}\left(  t_{2}\right)  ,f_{n}\left(
t_{1}\right)  \right)  \leq\left\vert t_{1}-t_{2}\right\vert \cdot
\ell(\mathfrak{C}_{n})\leq\left\vert t_{1}-t_{2}\right\vert \cdot M
\]
for all $t_{1},t_{2},n$, where $M$ is a constant. Hence one can apply the
Arzela-Ascoli theorem: there is a path $\mathfrak{C}$ given by
$f:[0,1]\rightarrow E_{2}$ such that $f_{n}\rightarrow f$ uniformly. Then
$\ell(\mathfrak{C})\leq m$, and $\mathfrak{C}$ meets every tangent. Indeed,
assume the contrary. If $\mathfrak{C}$ does not meet $t$, then it is at a
positive distance from $t$, and in view of the uniform convergence
$f_{n}\rightarrow f$, $\mathfrak{C}_{n}$ will be at a positive distance from
$t$ for $n\geq N$, which is impossible. Consequently, there is a minimum for
Theorem 2. As for Theorem 1, each $\mathfrak{C}_{n}$ starts at the center of
the circle, and hence so does $\mathfrak{C}$.

\bigskip

4. \textit{Notation}. If $A,B,C,D$ are points of the plane, then $AB$ will
denote the straight segment between $A$ and $B$, $\left\vert AB\right\vert $
the distance between $A$ and $B$, and $ABCD\ldots$ the path $AB\cup BC\cup
CD\cup\ldots$. If $A\neq B$, let $d(AB$) be the straight line that connects
$A$ and $B$, and $r(AB)$ the half-line of $d(AB)$ that starts at $A$ and
passes through $B$. If $A\neq B\neq C$ and $r(BA)$ is not the half-line
opposite to $r(BC)$, then $a(ABC)$ will be the closed convex angular region
bounded by $r(BA)$ and $r(BC)$, and $\measuredangle(ABC)$ the angular measure
of $a(ABC)$. Therefore, $0\leq\measuredangle(ABC)<\pi$ always.

The convex envelope of a set $S$ will be denoted by $k(S)$. In particular, if
$A,B,C,\ldots$ are points, $k(ABC\ldots)$ will be the smallest convex polygon
that contains $A,B,C,\ldots$; $k(ABC)$ will be the triangle with vertices
$A,B,C$.

The circle under consideration will be $K$, and its interior and exterior,
$\operatorname{int}(K)$ and $\operatorname{ext}(K)$.

If a path is given by $f:[a,b]\rightarrow E_{2}$ and $X=f(t),Y=f(s)$, we write
$X<Y$ if $t<s$.

\bigskip

5. First, let $\mathfrak{C}$ be an arbitrary path given by $c:[a,b]\rightarrow
E_{2}$. The interval $[a,b]$ is a union of $c^{-1}(K)$, $c^{-1}%
(\operatorname{int}(K))$ and $c^{-1}(\operatorname{ext}(K))$; as we know from
topology, $c^{-1}(\operatorname{int}(K))$ and $c^{-1}(\operatorname{ext}(K))$
are composed of open intervals in $[a,b]$, $c^{-1}(K)$ is composed of closed
intervals in $[a,b]$, and additionally, if there are infinitely many such
intervals, of the accumulation points. I will denote the paths corresponding
to the intervals of $c^{-1}(K)$, $c^{-1}(\operatorname{int}(K))$ and
$c^{-1}(\operatorname{ext}(K))$ by $j_{1},j_{2},j_{3}$ respectively, after
having assigned to them their initial and terminal points if necessary.

\bigskip

6. Consider now a $j_{1}$ of a minimal curve $\mathfrak{C}$. A\ $j_{1}$ is an
arc of $K$, but not all of $K$. If the ends of the arc are not the initial and
terminal point of $j_{1}$, then a part of the arc is traversed twice, which
allows for a shortening by a chord. Therefore, for a minimal $\mathfrak{C}$
the $j_{1}$ are arcs of $K$ traversed once. The $j_{2}$, being in the interior
of the circle, where there are no tangents, are obviously chords or segments
of chords.

\bigskip

7. It remains to consider the $j_{3}$. First let $X\in\operatorname{ext}(K)$;
let $t_{+}$ and $t_{-}$ be the two tangents of $K$ from $X$ such that, when
viewed from $X$, $K$ is on the left of $t_{+}$. Let $P_{+}(X)$ and $P_{-}(X)$
be the points at which $t_{+}$ and $t_{-}$ touch $K$ (Figure 4). \ For $X\in
K$, we write $P_{+}(X)=P_{-}(X)=X$. Thus $P_{+}$ and $P_{-}$ are continuous
maps of $K\cup\operatorname{ext}(K)$ onto $K$. Each $P_{+}(j_{3})$ and
$P_{-}(j_{3})$ is a closed arc on $K$, and if $j_{3}$ has a point on $K$, then
$P_{+}(j_{3})\cup P_{-}(j_{3})$ is a closed arc on $K$. Now let $X_{0}%
\in\mathfrak{C}\cap\operatorname{ext}(K)$ and let $t$ be a tangent that
properly separates $K$ from $X_{0}$ so that $X_{0}\notin t$ . Suppose for
simplicity that $X_{0}$ is not initial or terminal on $\mathfrak{C}$. There
are $Y,Z$ on $\mathfrak{C}$, $Y<X_{0}<Z$, such that for $Y\leq X\leq Z$, $X$
is also properly separated from $K$ by $t$. Let $X_{1},X_{2},X_{3},X_{4}$ be
four points, $Y\leq X_{1}\leq X_{2}\leq X_{3}\leq X_{4}\leq Z$, for which the
minimum and maximum of $P_{+}$ and $P_{-}$ on the path $\left\{
X\in\mathfrak{C}:Y\leq X\leq Z\right\}  $ are attained. Then the polygonal
path $YX_{1}X_{2}X_{3}X_{4}Z$ meets the same tangents and is strictly shorter
than any other path that contains $Y,X_{1},X_{2},X_{3},X_{4},Z$ in the same
order. It follows that each $j_{3}$ is a succession of segments. We call $AB$
a maximal segment of $\mathfrak{C}$ if $AB$ is not part of a segment
$A^{\prime}B^{\prime}$ of $\mathfrak{C}$ that properly contains $AB$.%
\begin{figure}[ptb]%
\centering
\includegraphics[
height=3.2287in,
width=3.4429in
]%
{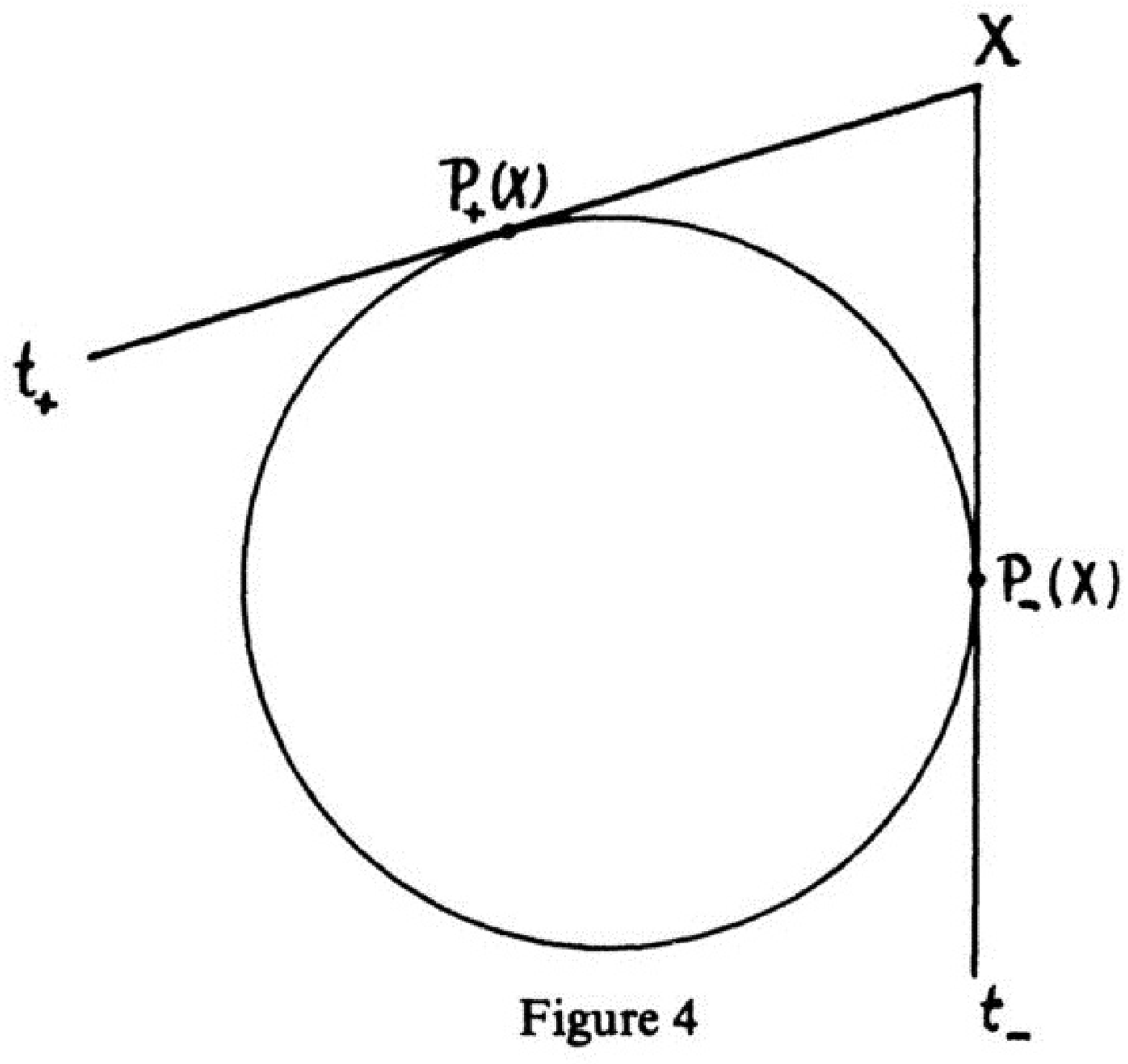}%
\end{figure}

\bigskip

8. It is useful to consider the problem in a different way. Let $D$ be a
convex compact domain on the plane $E_{2}$. The support lines of $D$ are the
lines $\ell$ that intersect $D$ in such a way that $D$ is contained completely
in one of the two half-planes defined by $\ell$. Let $S$ be a connected set.
Then $D$ belongs to $k(S)$ if and only if $S$ intersects all the support lines
of $D$. In particular, if $D$ is bounded by a smooth curve, then the support
lines are the tangents of the boundary curve.

Thus, our problem reduces to finding the shortest curve $\mathfrak{C}$ such
that $K\subseteq k(\mathfrak{C)}$.

We recall that $k(S_{1}\cup S_{2})$ is a union of segments $XY$ with $X\in
k(S_{1})$ and $Y\in k(S_{2})$.

We say that an arc $\mathfrak{S}$ of a path $\mathfrak{C}$ is repeating if
either $K\subseteq k(\mathfrak{C\smallsetminus S})$ or
$\mathfrak{C\smallsetminus S}$ intersects all the tangents to $K$. In
particular, $\mathfrak{S}$ is repeating if $\mathfrak{S}\subseteq
k(\mathfrak{C\smallsetminus S})$. \ A repeating arc is always a straight
segment traversed once.

\bigskip

9. Let $P\in K\cap\mathfrak{C}$, let $t$ be a tangent at $P$, and let
$Q\in\mathfrak{C}$ be properly separated from $K$ by $t$, that is, $Q$ belongs
to the interior of the half-plane defined by $t$ that does not contain $K$.
Then $\mathfrak{C}$ is a straight segment near $P$. Indeed, if $\mathfrak{S}$
is an arc of $\mathfrak{C}$ that contains $P$ and is contained in the interior
of the triangle $k\left(  P_{+}(Q)QP_{-}(Q)\right)  $ (Figure 5), then there
must be points $X,Y$ in the dashed regions with $X,Y\in
k(\mathfrak{C\smallsetminus S})$. Therefore, $\mathfrak{S}\subseteq
k(YQX)\subseteq k(\mathfrak{C\smallsetminus S})$ hence $\mathfrak{S}$ is
repeating, and hence a segment.%
\begin{figure}[ptb]%
\centering
\includegraphics[
height=3.8265in,
width=6.5894in
]%
{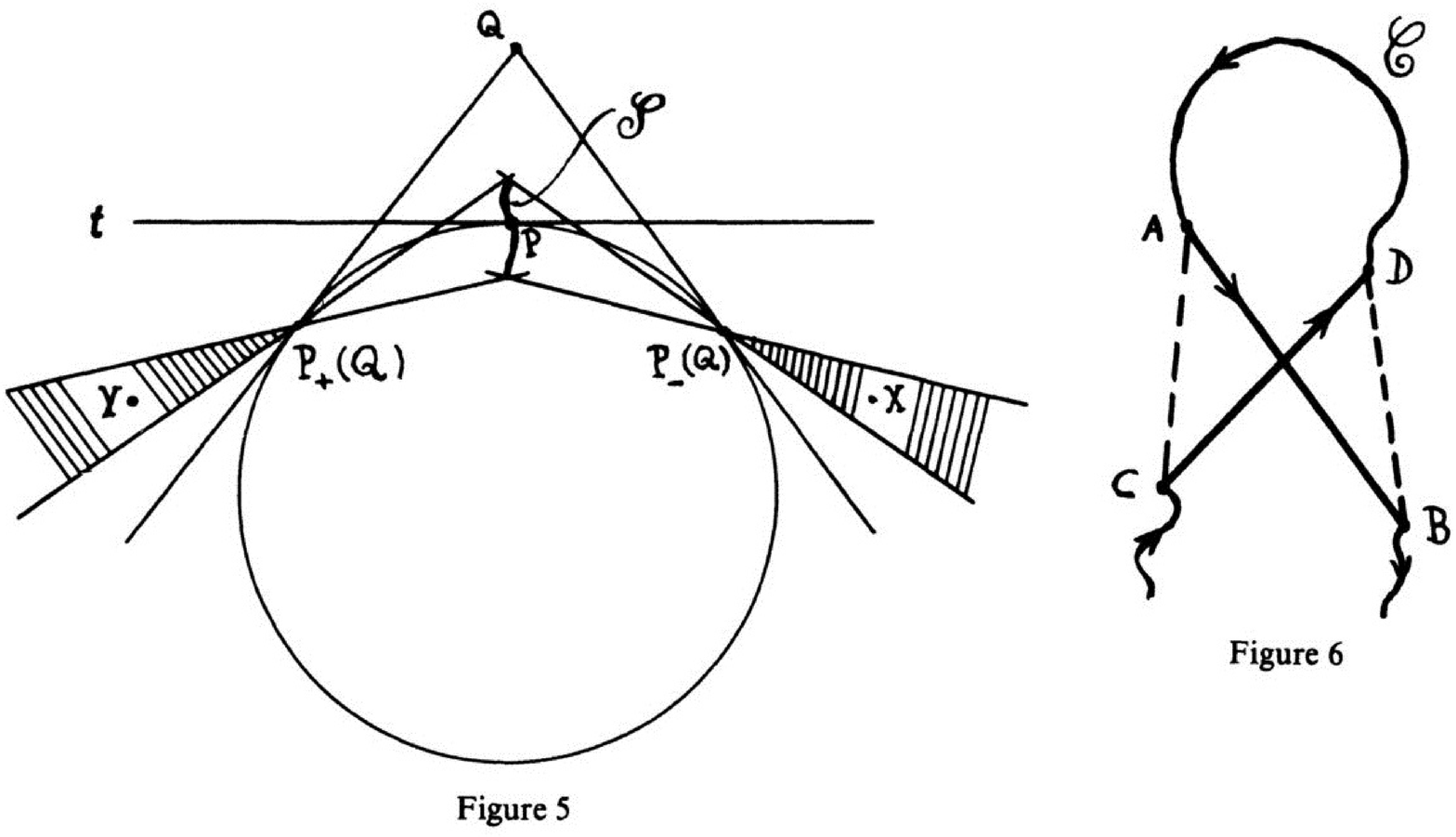}%
\end{figure}

\bigskip

10. The arguments of \S 8 allow us to exclude the following impossible situations:

Imp. a): Two noncollinear segments $AB$ and $CD$ that intersect at an interior
point of at least one of the segments $AB$ and $CD$. Indeed, in Figure 6,
$k(AB\cup DC)=k(DB\cup AC)$. We obtain a shortening by replacing $AB$ and $CD$
with $CA$ and $DB$ and by changing the direction of a part of $\mathfrak{C}$.

Imp. b): Two consecutive segments $ABC$ form an angle such that the opposite
angular region contains a point $X$ of $K$ or of $k(\mathfrak{C}\smallsetminus
ABC)$, $X\neq B$. Indeed, in Figure 7, let $S=(\mathfrak{C}\smallsetminus
ABC)\cup A\cup C$, $\mathfrak{C}=ABC\cup S$. There exists $Q\in k(S)$ such
that $X\in PQ$, $P\in k(ABC)$; then $ABC\subseteq k(QAC)$, $ABC$ is repeating
and must be replaced with the segment $AC$, which is shorter than $ABC$.%
\begin{figure}[ptb]%
\centering
\includegraphics[
height=2.7215in,
width=2.0473in
]%
{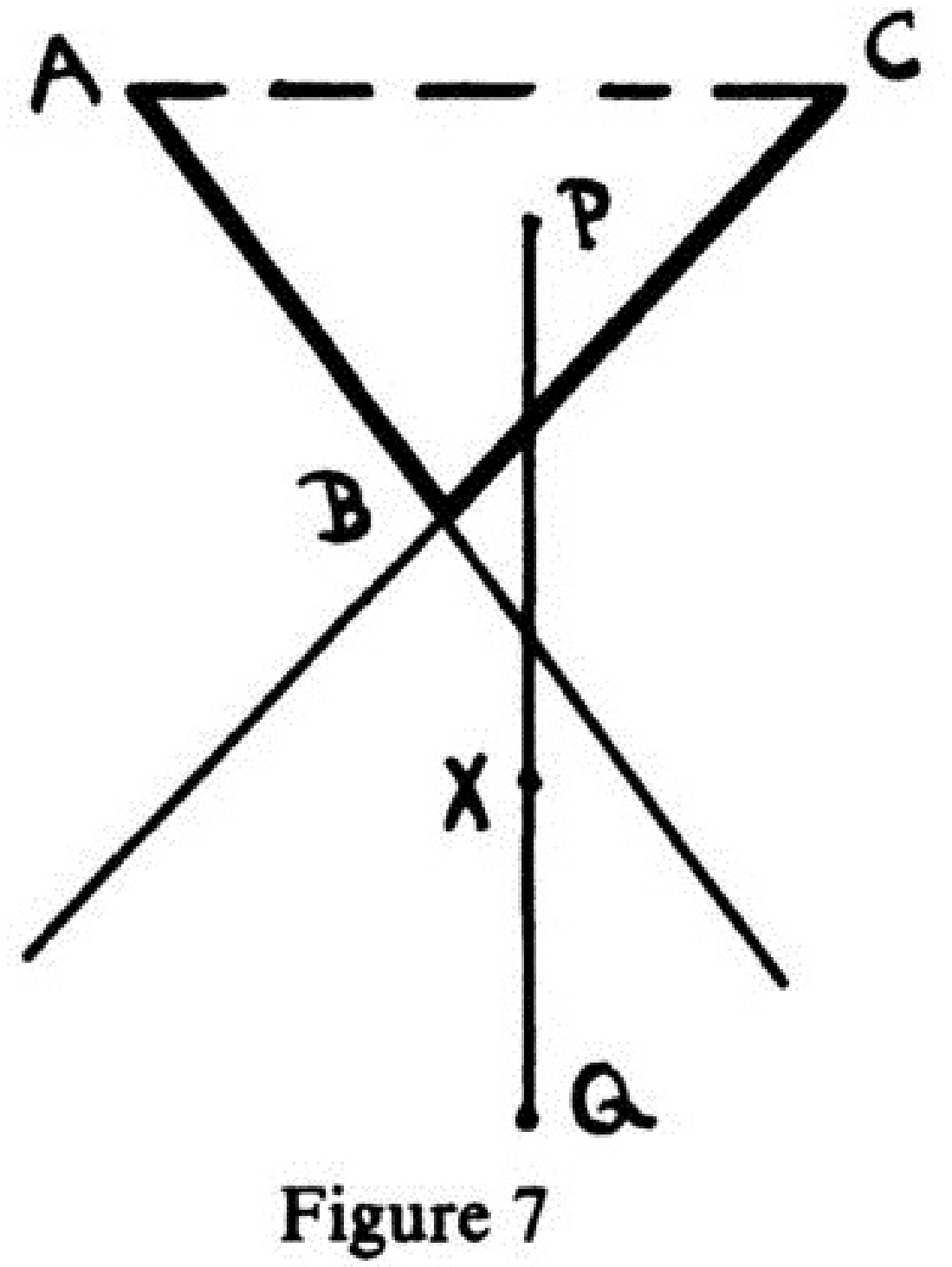}%
\end{figure}

Imp. c): A segment $AC$ such that $C$ is the initial or terminal point of
$\mathfrak{C}$ and such that $r(AC)$ intersects $K$ beyond $C$. The only
exception: $C$ is the initial point imposed by Problem (P).

Imp. d): Two coinciding segments at least one of which is followed by a
segment, in a different direction (Figure 8).%
\begin{figure}[ptb]%
\centering
\includegraphics[
height=3.0585in,
width=2.5745in
]%
{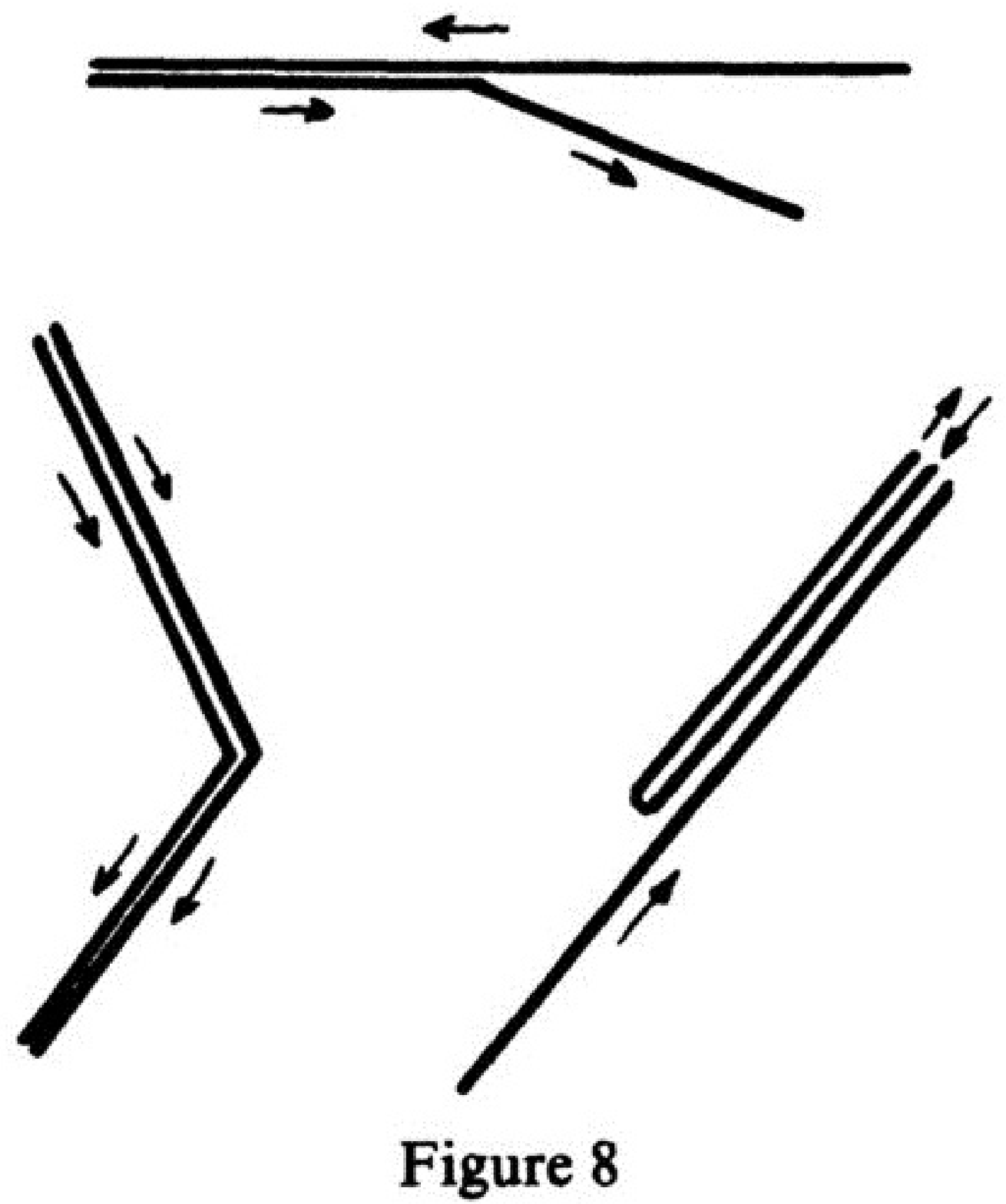}%
\end{figure}

Imp. e): Two consecutive segments $A^{\prime}BC^{\prime}$ such that $K$ is
contained in $a(A^{\prime}BC^{\prime})$. Indeed, if $A$ and $C$ are in the
interior of $A^{\prime}B$ and $C^{\prime}B$ respectively, sufficiently close
to $B$, then it is easily seen that $d(AC)$ properly separates $B$ from $K$
and that $ABC$ can be replaced with $AC$, which is shorter.

Imp. f): Two consecutive segments $ABC$ such that $K$ is contained in one of
the angular regions supplementary to $a(ABC)$, respectively, the half-planes
defined by $d(AB)$ if $r(BA)=r(BC)$. The only exception is the case that
$A=C\in K$ and $d(AB)$ is a tangent to $K$. We consider here the nondegenerate
case where $r(AB)\neq r(BC)$ and $ABC\subseteq\operatorname{ext}(K)$ (Figure
9). Suppose that $BA$ and $BC$ are maximal and that $d(BA)$ separates $K$ from
$BC$. Let the tangents $t_{1}$ and $t_{2}$ from $A$ and $B$ be as in the
figure. For $C$ we have the possibilities $C_{0},C_{1},C_{2}$. $BC_{2}$ is
repeating, hence $C\neq C_{0}$. If $C=C_{1}$, then $BA^{\prime}$ is repeating,
and hence $BA$ is not maximal. Therefore, $C=C_{2}$. If $C$ is the terminal
point of $\mathfrak{C}$, then $BC$ is superfluous (because the only imposed
terminal point of $\mathfrak{C}$ is the center of $K$). Consider the
continuation of $\mathfrak{C}$ beyond $C$. If this is a segment $CD_{1}$ with
$D_{1}$ properly separated from $K$ by $t_{1}$, then $BA$ cannot be maximal.
If, on the other hand, the continuation is $CD_{2}$ with $D_{2}$ in the same
closed half-plane (determined by $t_{1}$) as $K$, then $BCD_{2}$ is repeating,
which is impossible.%
\begin{figure}[ptb]%
\centering
\includegraphics[
height=3.2976in,
width=5.0278in
]%
{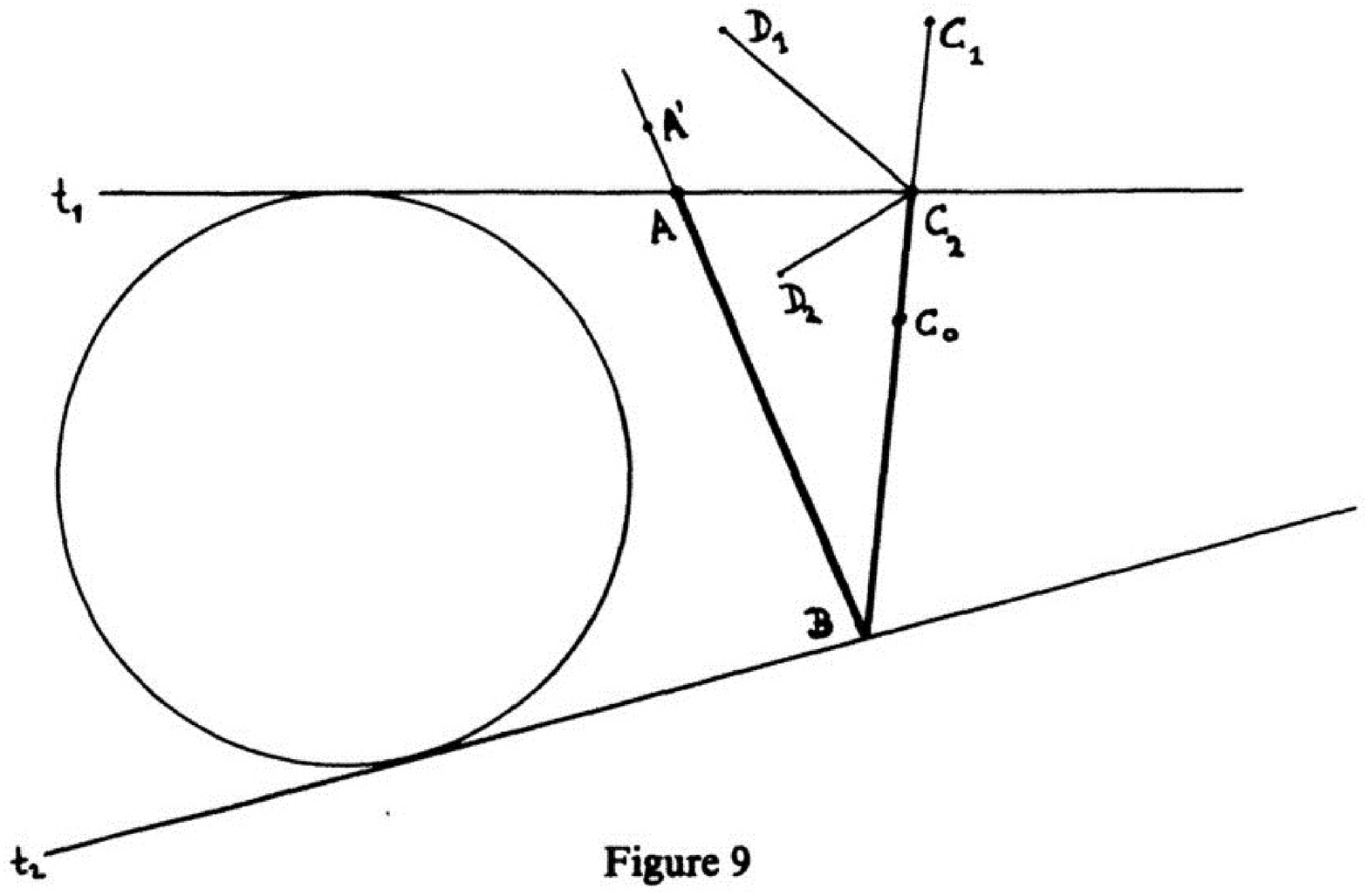}%
\end{figure}

\bigskip

11. We are now ready to determine what form $j_{3}$ can have. Let $AB$ be a
maximal segment of $j_{3}$ such that $d(AB)$ does not intersect $K$. In view
of Imp. b), e), and f), either $B$ is the terminal point or $\mathfrak{C}$
passes after $B$ along a segment directed towards $K$. The same is true for
$A$. If $AB\subseteq j_{3}$ is such that $d(AB)$ intersects $K$ [say, $r(AB)$
intersects $K$], then by Imp. c) and b), $B$ must be on $K$. If, moreover,
$d(AB)$ is not the tangent to $K$ from $B$, then the segment $AB$ must
continue to $\operatorname{int}(K)$ (in view of \S 9), beyond $B$. $A$ can be
terminal. If $A$ is not terminal, $\mathfrak{C}$ continues beyond $A$ along
$AC$, $C\in K$, or along $AC$ such that $d(AC)\subseteq\operatorname{ext}(K)$,
which is the case considered above.

Thus, we found for $j_{3}$ the seven possibilities indicated in Figure 10.%
\begin{figure}[ptb]%
\centering
\includegraphics[
height=5.1332in,
width=6.2233in
]%
{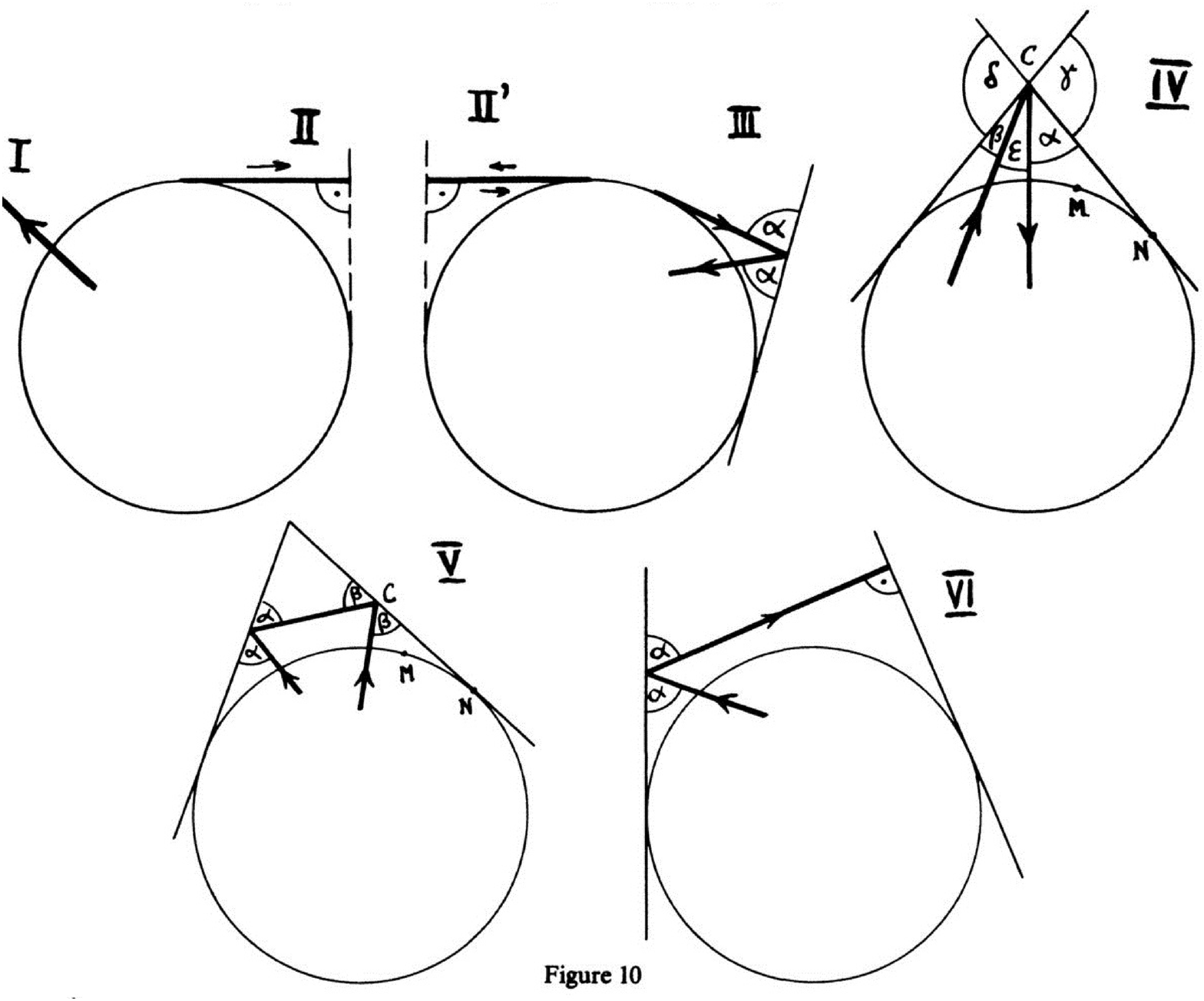}%
\end{figure}

The continuations to the interior of $K$ are indicated. Identical angles are
marked by identical letters. Also the right angles are indicated.

For example, the equality of the two angles $\beta$ in case V is due to the
following elementary fact: if $X$, $Y$ are two points on the same side of a
straight line $d$, then the shortest path between $X$ and $Y$ that passes
through $d$ consists of two segments $XZ$, $YZ$ with $Z\in d$ that form the
same angle with $d$. (The explanation of the right angles is even simpler.) In
case IV we obtain $\varepsilon=0$ as a limit case. It is easily seen that%
\begin{equation}%
\begin{array}
[c]{ccc}%
\alpha\leq\beta+\gamma, &  & \beta\leq\alpha+\gamma.
\end{array}
\end{equation}
It can also be shown that $2\beta+2\alpha+\varepsilon\geq\pi$, but we will not
use this inequality.

\bigskip

12. We have shown that $j_{1}$, $j_{2}$, and $j_{3}$ are the simplest arcs of
$\mathfrak{C}$; it remains to show that complicated configurations cannot be
produced by accumulations of $j_{1}$, $j_{2}$, $j_{3}$ on $K\cap\mathfrak{C}$.
First we prove the following.

If $AB$ is a $j_{2}$ of $\mathfrak{C}$, that is, a chord or part of a chord of
$K$, and $B\in K$ , then $B$ is not the terminal point of $\mathfrak{C}$, and
$\mathfrak{C}$ continues to a segment $BC$; $C\in\operatorname{ext}(K)$, $B\in
AC$, that is, $\mathfrak{C}$ continues straight ahead to $d(AB)$. First of
all, if $B$ were terminal, let $A^{\prime}\in AB\cap\operatorname{int}(K)$,
$\mathfrak{S}=(\mathfrak{C}\smallsetminus A^{\prime}B)\cup A^{\prime}$.
$\mathfrak{S}$ meets all the tangents, except perhaps those from $B$. But
since $\mathfrak{S}$ is closed, it meets also those from $B$ by continuity,
and hence $A^{\prime}B$ would be superfluous. In the same way, the rest of
$\mathfrak{C}$ is not formed by another $j_{2}$. The arguments that follow
refer to Figure 11: $t$ is the tangent from $B$ that we take to be horizontal.
If there is a $Y\in\mathfrak{C}$ that is properly separated from $K$ by $t$,
then the assertion follows from \S 9. Otherwise there will be a $Z>B$ such
that $\mathfrak{S=}\left\{  X\in\mathfrak{C}:Z\geq X\geq B\right\}  $ lies in
the rectangle $k(PQRW)$, where $R,W\in K$ and $P,Q\in t$, with small
$\left\vert QB\right\vert $ and $\left\vert PB\right\vert $. Suppose
$\measuredangle(A^{\prime}BP)\leq\pi/2$. To pass from $k(QBTR)$ to $k(PBTW)$,
$\mathfrak{S}$ must pass through $B$, since $\mathfrak{S}$ cannot cross
$A^{\prime}B$ inside, in view of Imp. a). If $\mathfrak{S}$ passes the two
sides of $A^{\prime}B$ infinitely many times, it must have infinitely many
loops $b$ departing from and arriving at $B$, in $k(PBTW)$. \ Let $b$ be such
a loop, with horizontal elongation $\left\vert BS\right\vert $. For a
sufficiently small $BP$, the path $A^{\prime}SB$ is shorter than $A^{\prime}B$
and $b$ together, since%
\begin{align*}
\left\vert A^{\prime}S\right\vert +\left\vert BS\right\vert  &  \leq
\frac{\left\vert A^{\prime}B\right\vert }{\cos\gamma}+\left\vert BS\right\vert
\leq\frac{\left\vert A^{\prime}B\right\vert }{\cos\gamma}-\left\vert
BS\right\vert +\operatorname{length}(b)\\
&  \leq\frac{\left\vert A^{\prime}B\right\vert }{\cos\gamma}-\left\vert
A^{\prime}B\right\vert \sin\gamma+\operatorname{length}(b)<\left\vert
A^{\prime}B\right\vert +\operatorname{length}(b)
\end{align*}
if $\gamma$ is sufficiently small. Therefore, if $b$ is a loop with maximum
horizontal elongation, then we can replace $A^{\prime}B$ with $A^{\prime}SB$
and omit all the loops in $k(PBTW)$. Therefore, we may assume that
$\mathfrak{S}$ is entirely on one side of $A^{\prime}B$, say in $k(QBTR)$.
However, in this case $\mathfrak{S}\subseteq k(QBR)$. Indeed, if $X$ is a
point in the triangle $k(BTR)$ without $RB$, then $X$ is on a chord $j_{2}$
that must intersect $RW$ or $BT$ in the interior, unless $\mathfrak{S}%
\subseteq j_{2}$, the case which has already been excluded. Now let
$U\in\mathfrak{S}$ be the leftmost for all of $\mathfrak{S}$. We replace
$A^{\prime}B\cup\left\{  X\in\mathfrak{S}:U\geq X\geq B\right\}  $ with
$A^{\prime}VU$, which meets all the tangents met by $A^{\prime}B\cup\left\{
X\in\mathfrak{S}:U\geq X\geq B\right\}  $ and has the length $\geq\left\vert
A^{\prime}B\right\vert +\left\vert BV\right\vert $, whereas for $\omega
=\measuredangle(A^{\prime}BV)$,%
\begin{align*}
\left\vert AA^{\prime}V\right\vert +\left\vert VU\right\vert  &  =\left\vert
A^{\prime}B\right\vert \cos\delta+\left\vert VB\right\vert \cos(\pi
-\omega-\delta)+\left\vert VU\right\vert \\
&  \leq\left\vert A^{\prime}B\right\vert \cos\delta+\left\vert VB\right\vert
\left\vert \cos(\omega+\delta)\right\vert +\left\vert VB\right\vert \tan
\alpha<\left\vert A^{\prime}B\right\vert +\left\vert VB\right\vert
\end{align*}
if $\alpha$ and $\delta$ are sufficiently small. This proves the
assertion.\bigskip%
\begin{figure}[ptb]%
\centering
\includegraphics[
height=3.1025in,
width=5.4545in
]%
{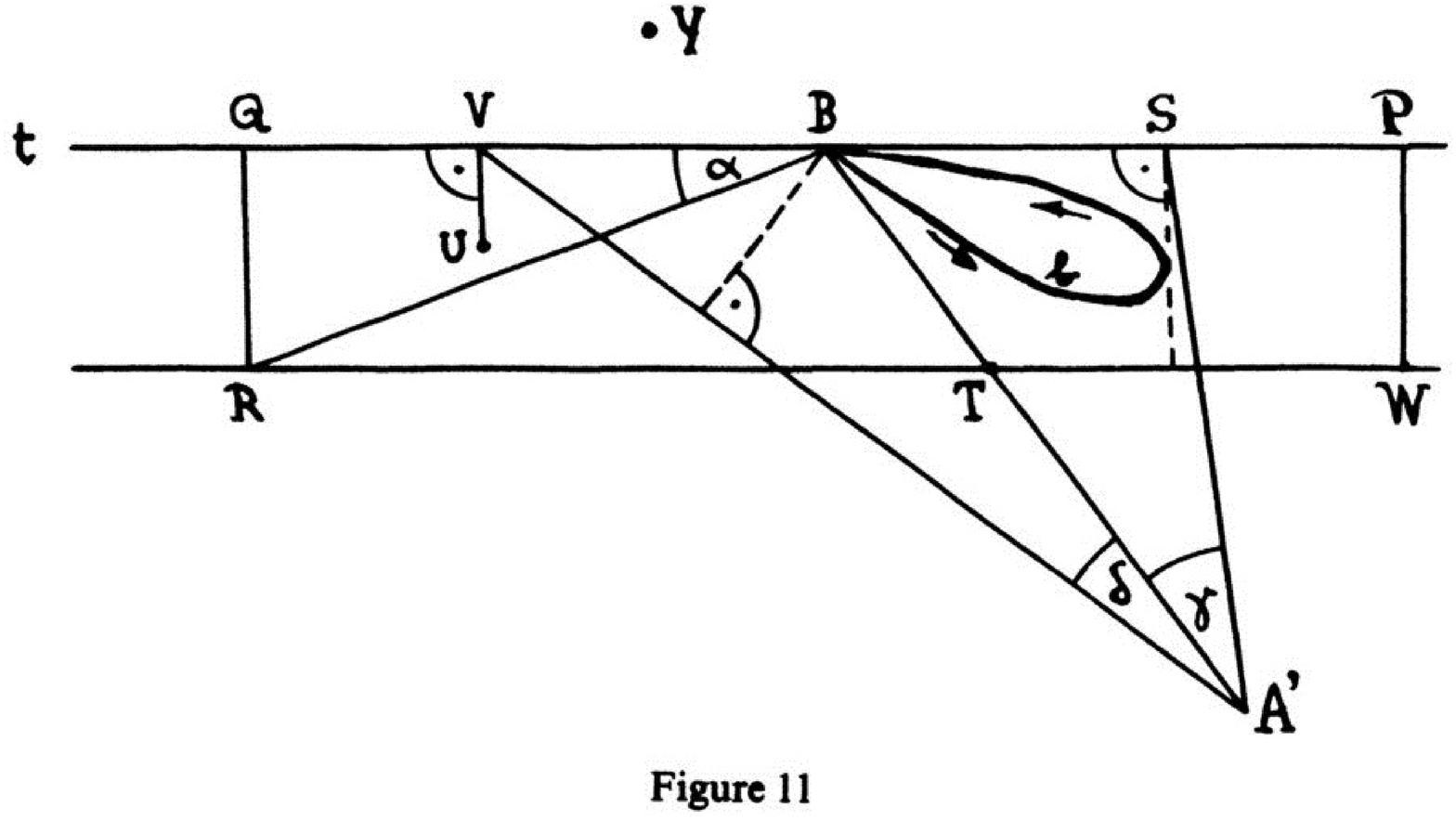}%
\end{figure}

13. Each $j=j_{1}$ or $j_{3}$ \textquotedblleft covers\textquotedblright\ a
set of tangents whose intersection points with $K$ form an arc $P_{+}(j)\cup
P_{-}(j)$ on the circle $K$. We denote this arc $P(j)$. We claim that $P(j)$
does not have common interior points with $P(j^{\prime})$ if $j\neq j^{\prime
}$. This is clear if $j$ or $j^{\prime}$ is a $j_{1}$, and hence an arc of
$K$, since a part of this $j_{1}$ can be replaced by a shorter chord. Thus,
assume that $j$ and $j^{\prime}$ are $j_{3}$. Obviously, $P(j)\varsubsetneq
P(j^{\prime})$; otherwise $j$ is repeating. Therefore, we are in a situation
of Figure 12. If $j$ is a $j_{3}$ of type V (Figure 10), we can shorten the
path by cutting through the angle formed by $\mathfrak{C}$ at $C$. If $j_{3}$
is of type IV (Figure 10), we can go down from $C$ along the left tangent,
which will ensure a shortening unless $\beta=\alpha+\gamma$, and so forth. We
arrive at the following situation as the only possibility (Figure 13): $j$
contains the segment $b=BD$, with $\beta\leq\pi/2$, and $j^{\prime}$ contains
the segment $a=AC$ with $\alpha\leq\pi/2$. In the same way as in the proofs of
\S 10, if $S=\mathfrak{C}\smallsetminus a\smallsetminus b$, we find points
$X,Y,W\in k(S)$, in the indicated angular regions for which $R,Q,P$,
respectively, are in $k(\mathfrak{C})$. But then $a,b\subseteq k(ABWXY)$, and
$a\smallsetminus A$ and $b\smallsetminus B$ are repeating, as well as their
continuations up to the tangents $t_{1}$ and $t_{2}$, respectively, which
would enforce a forbidden crossing. We note that this is the first time we
used the fact that $K$ is a circle, or rather that the normals to the tangents
at $Q$ and $P$ meet at $Z\in\operatorname{int}(K)$. Up until this moment,
everything was applicable for smooth and convex $K$s that do not contain
straight segments.\bigskip%
\begin{figure}[ptb]%
\centering
\includegraphics[
height=1.8738in,
width=3.0444in
]%
{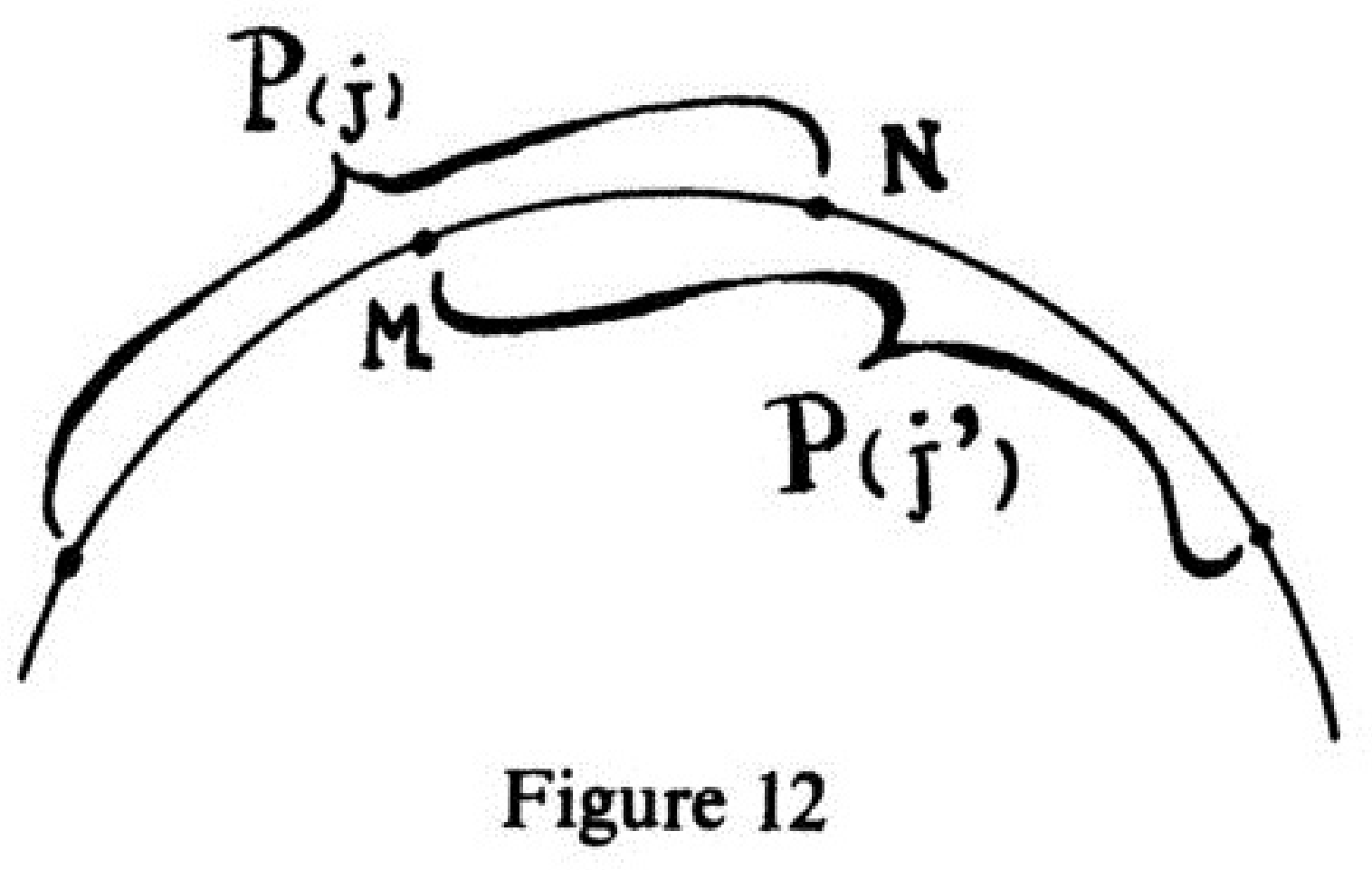}%
\end{figure}
\begin{figure}[ptb]%
\centering
\includegraphics[
height=4.7638in,
width=6.0216in
]%
{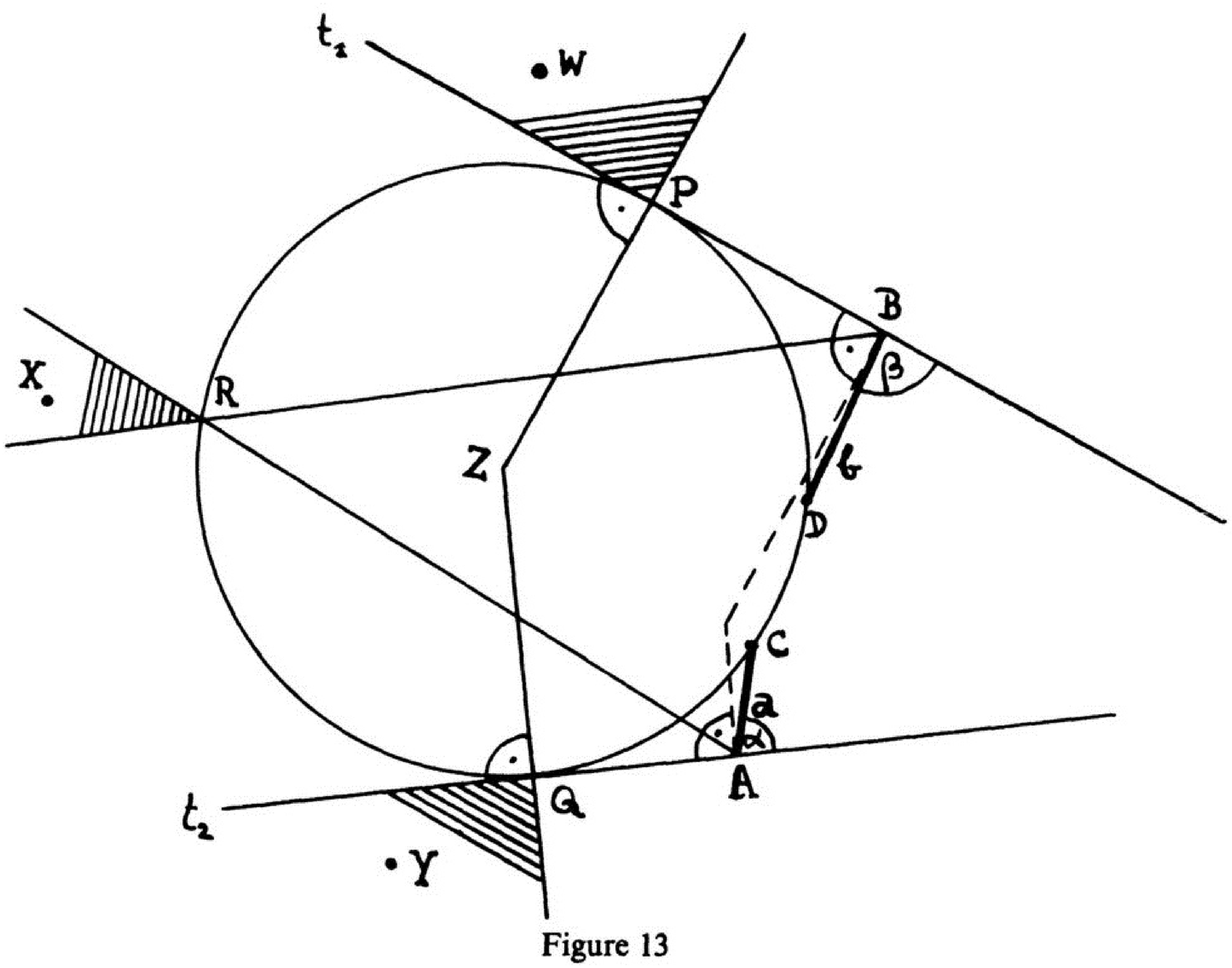}%
\end{figure}

14. It is easily seen from the above that if $j=j_{1}$ or $j_{3}$, then
$\mathfrak{C}\smallsetminus j$ is entirely on one side of $K$ with respect to
each tangent that touches $K$ at $P(j)$, that is, in the non-dashed part of
the plane shown in Figure 14. It follows that if $AB\subseteq\mathfrak{C}$ is
a maximal segment that contains a chord $PQ$ of $K$, then neither $A$ nor $B$
is terminal and the continuations of $\mathfrak{C}$ from $A$ and from $B$ do
not go in the same half-plane determined by $d(AB)$. In other words, we have
the situation of Figure 15.%
\begin{figure}[ptb]%
\centering
\includegraphics[
height=3.2179in,
width=6.6517in
]%
{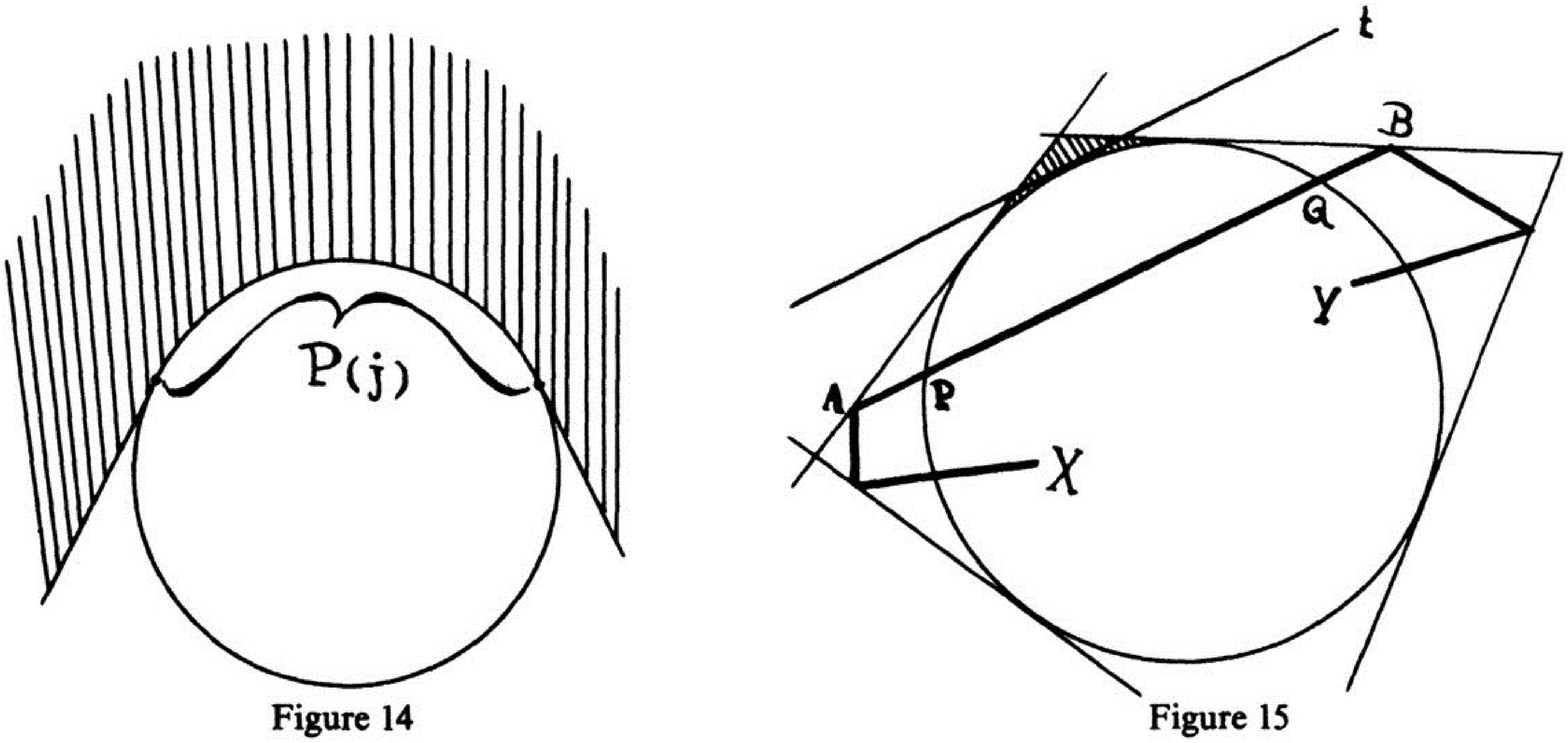}%
\end{figure}

The intersection of $\mathfrak{C}$ with the tangent $t\parallel AB$ must be in
the dashed region, which $\mathfrak{C}$ can reach neither from $X$ nor from
$Y$ without crossing $AB$.\bigskip

15. We show that $\mathfrak{C}$ consists of finitely many $j_{1},j_{2},j_{3}$.
Assume the contrary. Then there exists an $X\in K\cap\mathfrak{C}$ such that
for all $Y,Z\in\mathfrak{C}$, $Y<X<Z$, there are infinitely many $j_{1}%
,j_{2},j_{3}$ between $Y$ and $Z$, say between $Y$ and $X$. Suppose there are
no chords $j_{2}$ among them. Then the $j_{3}$ must be of type II' in Figure
10. But each of those has length $2$, and so there are only finitely many of
them; therefore, if $Y$ is close enough to $X$, there are only $j_{1}$'s,
whence $\left\{  W\in\mathfrak{C}:Y\leq W\leq X\right\}  $ is an arc on $K$
and belongs to a single $j_{1}$. Consequently, there are infinitely many
chords $j_{2}$ approaching $X$. In view of \S 14, they must zigzag, as in
Figure 16. The $j_{2}$ have exactly the same direction as the tangent $t$ at
$X$. To each $j_{2}$ we attach a $j_{3}$ of the form III, IV, or V (Figure
10).%
\begin{figure}[ptb]%
\centering
\includegraphics[
height=3.8265in,
width=6.2349in
]%
{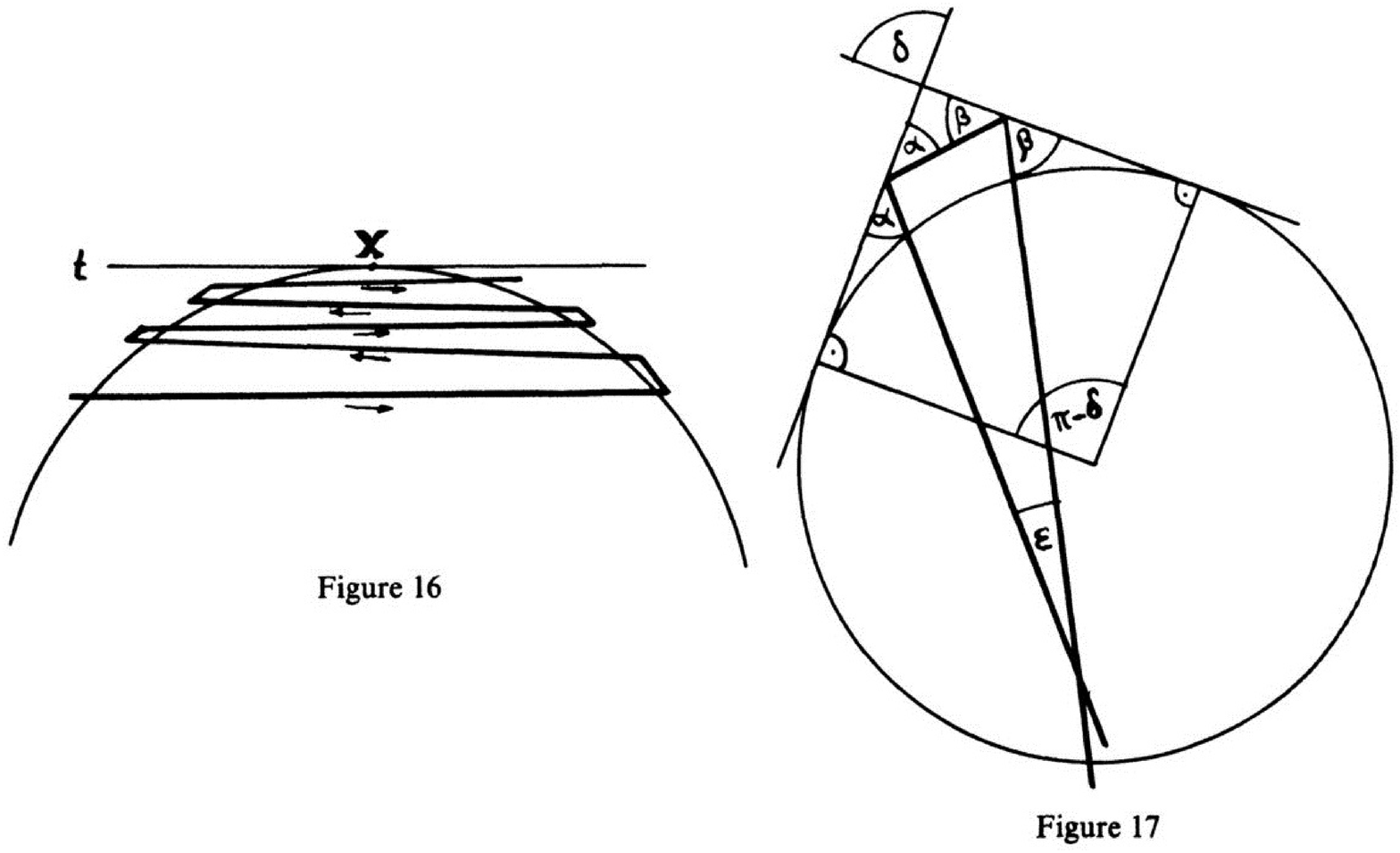}%
\end{figure}

For a $j_{3}$ of the form V we have, in Figure 17, $2\pi=\varepsilon
+\delta+\pi-\alpha+\pi-\delta$; therefore $\alpha+\beta=\varepsilon+\delta$.
But $\alpha+\beta=\pi-\delta$, and hence $\delta=\pi/2-\varepsilon/2$; for
$\left\vert \varepsilon\right\vert $ small, $\delta\approx\pi/2$, and
$P(j_{3})$ has a length $\approx\pi/2$, which is too large. In the same way we
have contradictions for $j_{3}$ of the types III and IV. Thus we have found
that $\mathfrak{C}$ is a finite union of paths $j_{1},j_{2},j_{3}$.\bigskip

16. We can now prove Theorems 1 and 2. In Theorem 1 we have one free terminal
point; in Theorem 2, both terminal points are free. We see that an arc $j_{1}$
on $K$ cannot be a terminal path of $\mathfrak{C}$. We already know that the
$j_{2}$ are not terminal except for the fixed terminal point in Theorem 1. In
view of \S 14, the $j_{3}$ of type I are not terminal. Therefore, only the
$j_{3}$ of type II or VI remain candidates for a free terminal path.

Suppose that the free terminal segment is a $j_{3}$ of type II, and hence a
tangent of length $1$. It may be followed by a $j_{3}$ of type III or by an
arc $j_{1}$. In the first case we have Figure 18. But here we see that $ABC$
can be replaced with $BAC$, which is shorter and has the same convex envelope.
The same argument is valid if $\mathfrak{C}$ ends with a $j_{3}$ of type VI.
Consequently, the only possibility is a tangent of length $1$ followed by an
arc. We show that $\mathfrak{C}$ does not admit entire chords. Indeed, assume
the contrary and take the first such chord after the arc. We obtain Figure
19.
\begin{figure}[ptb]%
\centering
\includegraphics[
height=2.7754in,
width=6.2748in
]%
{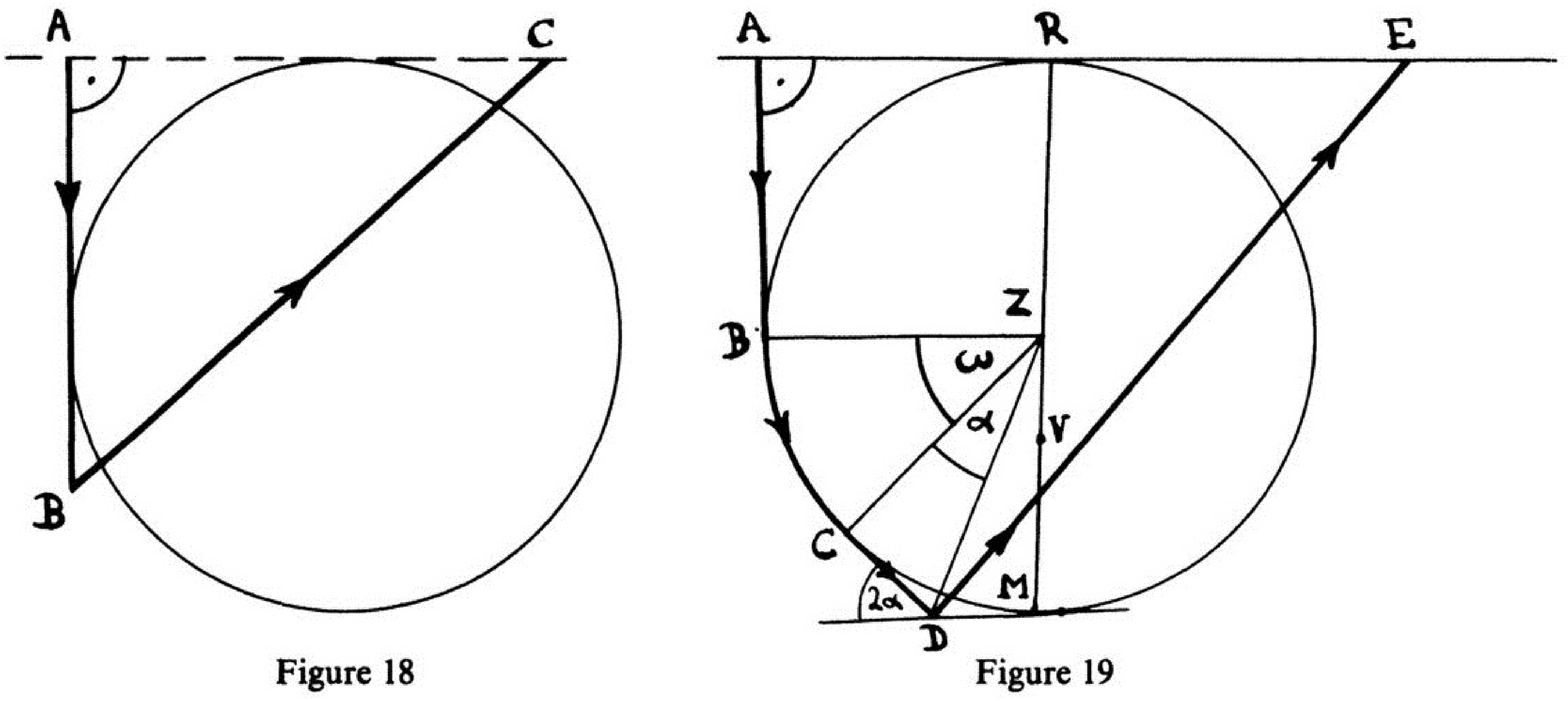}%
\end{figure}
Let $\alpha+\omega\leq\pi/2$. If%
\[
\left\vert RE\right\vert <\frac{4}{\pi-2}-\frac{\pi-2}{4},
\]
we verify that the circle centered at $E$ of radius $\left\vert ER\right\vert
+(\pi/2)-1$ passes through $V$ (on the diameter $RM$), and hence%
\begin{equation}
\left\vert ED\right\vert >\left\vert ER\right\vert +\frac{\pi}{2}-1.
\end{equation}
This enables us to replace the path $ABCDE$ (where $BC$ is an arc) with
$DCBRE$, where $CBR$ is an arc on $K$. The length of the first is
$1+\omega+\tan\alpha+\left\vert DE\right\vert >1+\omega+\tan\alpha+\left\vert
ER\right\vert +(\pi/2)-1=\omega+\tan\alpha+(\pi/2)+\left\vert ER\right\vert $
which is the length of the second. If%
\[
\left\vert RE\right\vert >\frac{4}{\pi-2}-\frac{\pi-2}{4},
\]
we consider Theorem 2 first. The length of the path will be at least%
\[
\left\vert AB\right\vert +\left\vert BE\right\vert >2+\left\vert RE\right\vert
\geq2+\frac{4}{\pi-2}-\frac{\pi-2}{4}>\pi+2,
\]
and hence it is longer than the path in the theorem.

For Theorem 1, the curve must return to the center from $E$. Therefore, the
path will be longer than $\left\vert AB\right\vert +\left\vert BE\right\vert
+\left\vert EZ\right\vert >2\left(  1+\left\vert RE\right\vert \right)
>2+2\pi$, and hence longer than the path of Theorem 1.

Now let $\omega+\alpha>\pi/2$. For Theorem 1, we see that $DE$ is an obstacle
on the way of $\mathfrak{C}$ back to the center. For Theorem 2, we start from
the other terminal point, but we must have $\omega^{\prime}+\alpha^{\prime
}<\pi/2$ in view of \S 14, which reduces to the previous case.

Thus, we see that the single $j_{2}$ can only be the semichord leading to the
center, for Theorem 1, and so we cannot obtain anything different from the
paths described in the theorems.\bigskip

\textbf{Remarks}

(a) \ The problem considered can be largely generalized by replacing, for
example, $K$ with an arbitrary convex domain, or the tangent lines with
tangent circles, etc. The general statement of the problem is this: In an
arcwise connected metric space, given a family $\mathfrak{F}$ of closed $F$
and a connected compact set $K$ that intersects each $F\in\mathfrak{F}$, find
the shortest path that intersects all the $F$.

(b) \ In higher dimensions, if $E_{n}$ is $n$-dimensional Euclidean space,
$S_{n-1}$ the unit sphere, and $\ell_{n}$ the length of the shortest path that
meets all hyperplanes tangent to $S_{n-1}$, or the shortest path
$\mathfrak{C}$ with $S_{n-1}\subseteq k(\mathfrak{C)}$, then by induction%
\[
\ell_{3}\geq\sqrt{\left(  2+\sqrt{3}+\frac{7}{6}\pi\right)  ^{2}+4}%
\approx7.6628,
\]%
\[
\ell_{n}\geq\operatorname{const}+2n.
\]
By constructing particular paths we find that%
\[
\ell_{n}\leq\operatorname{const}\cdot\,n^{3/2},
\]%
\[
\ell_{3}\leq4+\frac{1}{2}\sqrt{2}\cdot3\cdot\pi\approx10.6643.
\]
The upper bounds seem to me to be the closest to reality.

(c) \ In an infinite-dimensional Hilbert space, say $H=l^{2}$, the
circumstances are slightly different. If $S$ is the unit sphere, there is no
path of finite length, nor a compact path $\mathfrak{C}$ such that $S\subseteq
k(\mathfrak{C)}$, since $k(\mathfrak{C)}$ is compact, but $S$ is not.
Therefore, we must consider compact convex sets $K\subseteq H$. In this case
there surely exists a compact path $\mathfrak{C}$ (which can be easily
constructed) such that $k(\mathfrak{C)}\supseteq K$. [In view of a theorem by
Hahn and Mazurkiewicz, there even exists a continuous map $f:[0,1]\rightarrow
H$ with $\operatorname{Im}f=K$ (see \cite{Newman}).] However, this curve in
general is not of finite length. If%
\[
K=\left\{  (x_{1},x_{2},\ldots):%
{\displaystyle\sum\limits_{j=1}^{\infty}}
j\,x_{j}^{2}\leq1\right\}  ,
\]
then $K$ is compact and convex, but any curve $\mathfrak{C}$ with
$k(\mathfrak{C)}\supseteq K$ has infinite length. If%
\[
K=\left\{  (x_{1},x_{2},\ldots):%
{\displaystyle\sum\limits_{j=1}^{\infty}}
j^{100}\,x_{j}^{2}\leq1\right\}  ,
\]
then there exists $\mathfrak{C}$ of finite length with $k(\mathfrak{C)}%
\supseteq K$. Moreover, it can be shown that if $K$ is compact and convex and
$\mathfrak{C}$ has finite length and $k(\mathfrak{C)}\supseteq K$, then there
exists a minimal curve. This can be done as in \S 3. It must simply be shown
that if $\left\{  \mathfrak{C}_{n}\right\}  $ is a minimal chain of curves,
then $%
{\textstyle\bigcup\nolimits_{n=1}^{\infty}}
\mathfrak{C}_{n}$ is relatively compact.

\begin{flushright}
Written by Henri Joris, Gen\`{e}ve

Translated by Natalya Pluzhnikov
\end{flushright}

\textbf{Addendum}

The only citation to the literature given in Joris (1980) is \cite{Newman}. \ 

I have included \cite{Be, Is, Me, Eg, FP, FM, Fi, FW, Gh} for the sake of
completeness. \ 

John Wetzel provided valuable comments as I\ prepared this draft for posting.

\end{document}